\documentclass[12pt]{article}
\usepackage[T1]{fontenc}
\usepackage{fixltx2e}
\usepackage{textcomp}
\usepackage{epsf}
\usepackage[dvips]{epsfig}
\usepackage{amssymb,amsmath,amsfonts,authblk}
\usepackage{srcltx}
\usepackage{graphicx}
\usepackage{epsfig}
\usepackage{graphics}
\usepackage{gastex}
\usepackage[usenames]{color}
\title{Words containing all permutations of a family of factors}
\author{Anna E. Frid}
\newtheorem{theorem}{Theorem}
\newtheorem{lemma}{Lemma}
\newtheorem{corollary}{Corollary}

\newtheorem{remark}{Remark}

\begin{document}
\maketitle
\begin{abstract} 
In this note, we prove that if a uniformly recurrent infinite word contains as a factor any finite permutation of words from an infinite family, then either this word is periodic, or its complexity (that is, the number of factors) grows faster than linearly. This result generalizes one of the lemmas of a recent paper by de Luca and Zamboni, where it was proved that such an infinite word cannot be Sturmian.
\end{abstract}

%\section{Introduction}
In a recent paper \cite{dlZ}, Aldo de Luca and Luca Q. Zamboni proved the following lemma:

\begin{lemma}[\cite{dlZ}] Let $r \in \mathbb N^+$ and $x \in \{0,1\}^{\omega}$ be a $r$-power-free Sturmian word. Then for each infinite sequence $V_0,V_1,\ldots$ with $V_i\in \{0,1\}^+$ there exist $k\geq 1$, $0\leq n_1 <n_2<\cdots<n_k$, and a permutation $\sigma$ of $\{1,2,\ldots,k\}$ such that $V_{n_{\sigma(1)}} V_{n_{\sigma(2)}}\cdots  V_{n_{\sigma(k)}} \notin \mbox{\rm Fact}(x)$. 
 \end{lemma}

This lemma was then used to prove the non-existence for any Sturmian word $x$ of a so-called $\varphi$-ultra monochromatic factorization, where $\varphi$ is a function on finite words equal to 0 if and only if the word is a factor of $x$ and 1 otherwise (for the definition of an ultra-monochromatic factorization, see again \cite{dlZ}).

In this short note, we extend the result by de Luca and Zamboni to all aperiodic uniformly recurrent words whose (factor) complexity is linear. The non-existence of a $\varphi$-ultra monochromatic factorization of such words follows then immediately, for the same $\varphi$ as before.

So, this is the result of this note.
\begin{theorem}
 Consider a uniformly recurrent infinite word $w$ and a family of its non-empty different factors 
$A=\{u_1,\ldots,u_n,\ldots\}$ such that $w$ contains as factor any permutation of a finite subset of words of $A$. Then either $w$ is periodic, or its factor complexity is greater than linear.
\end{theorem}
{\sc Proof.} 
For each subset $B$ of $A$ denote by $l_B$ the sum of lengths of words from $B$ (that is, the length of any permutation of them) and by $P_B$ the set of words equal to these permutations. So, $1\leq \#P_B \leq (\#B)!$. By our assertion, for any finite $B \subset A$, all words from $P_B$ are factors of $w$. The cardinality of $P_B$ is maximal if all permutations of words of $B$ are different, and minimal if all the words from $B$ commute, that is, if they all are integer powers of the same finite word.

Since the set $A$ is infinite, the length of its elements is unbounded. Since the word $w$ is uniformly recurrent, either it is (strictly) periodic, or the minimal period of words from $A$ is unbounded too. We suppose that the latter holds.

\begin{lemma}
For each finite subset $B$ of $A$ there exists a length $N$ such that the complexity $p_w(N)\geq N/4 \cdot (\#P_B-1)$.
\end{lemma}
{\sc Proof.}
Given a subset $B$ of $A$, choose another sufficiently long word $v$ from $A$: we need its length $|v|$ to be at least twice greater than $l_{B}$, and its minimal period $p_{min}$ to be greater than or equal to $l_{B}$. Clearly, such a word $v \in A$ exists. Let its length be $m$, and consider the set $F$ of all factors of $w$ of length $n=\lfloor 3/2 m\rfloor +l_{B}\leq 2m$. We shall give a lower bound to the number of these factors.

Indeed, consider only words from $F$ containing whole occurrences of words $vb$, where $b\in P_{B}$. Such a word is of length $m+l_{B}$ and thus should start at a position from 0 to $\lfloor m/2 \rfloor$ to be a factor of a word from $F$. For a given word $b \in P_{B}$ and a given starting position $k \in \{0,\ldots,\lfloor m/2 \rfloor\}$ of $vb$, consider a respective word $s(b,k)$ of length $n$. Note that in principle, there can be several words with these properties, but we just choose one of them. It exists since $vb$ is a factor of $w$ and $w$ is uniformly recurrent.

When is it possible that $s(b,k)=s(b',k')$ for some $b,b' \in P_B$ and $k,k'\in \{0,\ldots,\lfloor m/2 \rfloor\}$? If it is the case, the beginnings of (the fixed occurrence of) $v$ in this word are at the distance $p=|k-k'|\leq \lfloor m/2 \rfloor\}$ from each other. If $p=0$, then $k=k'$ and thus $b=b'$. Otherwise, suppose that $k<k'$. The word $v$ is $p$-periodic, and its minimal period $p_{min}$ is a divisor of $p$. So, $|b|=l_{B}\leq p_{min}\leq p$, and the word $b$ is equal to the prefix of length $l_{B}$ of the suffix of length $p_{min}$ of $v$. It means that $b$ (but not $b'$) is uniquely determined by $v$. So, the equality $s(b,k)=s(b',k')$ is possible only for one of $\#P_B$ possible words $b$. For all other elements $b'$ of $P_B$, the words $s(b',k)$ are all different for all $k \in \{0,\ldots,\lfloor m/2 \rfloor\}$. This gives
\[p_w(2m+4)\geq p_w(2m) \geq p_w(n)\geq (m/2+1)(\# P_B-1),\]
that is,
\[p_w(N)\geq N/4 \cdot (\# P_B-1)\]
for $N=2m+4$. \hfill $\Box$

\begin{corollary}
 If the complexity $p_w(n)$ is linear, then there exists a constant $C$ such that $\# P_B \leq C$ for all finite subsets $B$ of $A$.
\end{corollary}

Now let us take a set $B$ such that $\# P_B$ is maximal possible. Clearly, for any superset $B'\supset B$, we have $\# P_{B'}\geq \#P_{B}$ and thus these cardinalities are equal. Let us take $B'$ to be equal to $B \cup \{u',u''\}$, where $u'$ and $u''$ are any two elements of $A \backslash B$. We have $\#(P_Bu'u'' \cup P_Bu''u')\leq\#P_B$, that is, in fact the cardinalities are equal, which is possible only if $u'u''=u''u'$. In other terms, the words $u'$ and $u''$ commute. This is possible if and only if they are integer powers of the same primitive word $r$. Since this is true for any pair $u',u'' \in A \backslash B$, we see that all the elements of $A$ except for a finite number have the same minimal period $|r|$, which contradicts to the fact that their minimal period is unbounded.

This contradiction means that in fact, a non-periodic uniformly recurrent infinite word containing all finite permutations of words from an infinite set $A$ cannot have linear complexity. \hfill $\Box$

\medskip
\begin{remark}
{\rm In fact, we use not the fact that $w$ is uniformly recurrent but a weaker condition that periods of factors in $A$ can be arbitrarily large. I believe that if it is not the case, the result still holds.}
\end{remark}


\begin{thebibliography}{1}
 \bibitem{dlZ}
Aldo de Luca, Luca Q. Zamboni, On some variations of coloring problems of infinite words, J. Combin. Theory Ser. A 137 (2016) 166--178.
\end{thebibliography}
\end{document}